\documentclass{article}
\usepackage{amssymb,amsmath,amsthm,graphicx}
\usepackage[all,color]{xy}

\def\qed{\hfill {\hbox{${\vcenter{\vbox{               %HOLLOW SQUARE
   \hrule height 0.4pt\hbox{\vrule width 0.4pt height 6pt
   \kern5pt\vrule width 0.4pt}\hrule height 0.4pt}}}$}}}

\theoremstyle{definition}

\date{}

\title{\Large \textbf{Richard A. Litherland: A Brief Biography}}

\author{
Alison Litherland\footnote{Email: alison.litherland@gmail.com}\and
Patrick Litherland\footnote{Email: patricklitherland@outlook.com} \and
Sam Nelson\footnote{Email: Sam.Nelson@cmc.edu. Partially supported by Simons Foundation collaboration grant 702597}
\and
Steven Wallace\footnote{Email: steven.wallace@mga.edu}}

\begin{document}
\maketitle

\begin{abstract}
Richard A. Litherland was born in 1953 in England. He received his PhD at 
Trinity College in Cambridge in 1979 and moved to the USA in 1983. He had 
a lengthy and distinguished career as a professor of mathematics and 
researcher of low-dimensional topology,
based primarily at Louisiana State University (LSU) in Baton Rouge until
his untimely passing in November 2022. In this paper we (Rick's family 
and students) recount some memories of Rick, his life and his mathematics.
\end{abstract}

\parbox{5.5in} {\textsc{Keywords:} R.A. Litherland, low-dimensional topology,
knot theory

\smallskip

\textsc{2020 MSC:} 01A99, 57K10}

%\section{\large\textbf{Introduction}}\label{I}

\section{\large\textbf{Life Events}}\label{LE} %author: Alison and Patrick

As a teenager, Richard Litherland (1953-2022) inadvertently summed up 
something of his own character in a postcard home from Scout camp: ``Today 
I went on my 2nd class hike \dots but we lost the way and failed. I have 
passed Kim’s Game and knots''.

In two sentences we have: a hiker, a self-critic, an absorbent mind and a 
topologist.

Rick was born in 1953 in St Helens, a glass-making, Rugby--playing town not 
far from Liverpool in north-west England. He grew up first there and later 
in Flixton near Urmston, an outpost of Manchester between the River Mersey 
and the Manchester Ship Canal, where interwar suburbia meets the mosses and 
farmland of Cheshire. He was educated at Manchester Grammar School and 
Trinity College, Cambridge.

Those who knew him in later life may be surprised to find out what an athlete 
Rick was in his youth: a fly--half and sometime wing three-quarter on his 
grammar school's Rugby team, a regular at the school's Grasmere fellwalking 
camp, a ``Trek'' participant, the leader of smaller--scale family treks in the 
Peaks and Lakes, a Sea Scout. Yes, he got lost once or twice, and a Sea Scout 
expedition he was on had to be escorted home by the Royal Navy across a 
fog--bound English Channel, but he dined out on those failures and was proud 
of his successes, especially of completing the Pennine Way. In America, 
however, he was to relinquish outdoor pursuits apart from a habit of walking 
everywhere in car--addicted Baton Rouge.

In the headquarters of the Nelson's Sea Scouts troop, a pipesmoke-filled hut 
in the hinterland of Urmston where Nelson-era seafaring skills were preserved, 
Rick discovered knotting. It was in this period that he bought the copy of 
Clifford Ashley's topic--defining Book of Knots that would stay with him to 
the end of his life; the family home started to fill with Turk's heads and 
monkey's paws, and girlfriends began to receive gifts of macram\'{e} neckwear. 
Not just the look, but the structure of the knots fascinated him: the up and 
down, the over and under, the regularities and the variations.

\[\begin{array}{c}
\includegraphics{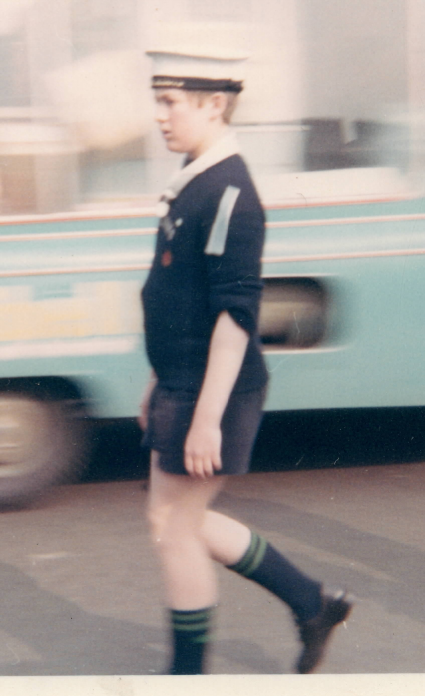}\\
\mathrm{R.\ A.\ Litherland\ on\ Scout\ parade,\ 1960s}
\end{array}\]

%[icecream-seascout-closer.jpeg could go here. Caption e.g.: R. A. Litherland on Scout parade, 1960s]

Meanwhile, the practical skill of knot tying was on a convergent trajectory 
with his love of mathematics. Although he did well at all school subjects 
apart from music and divinity, and although his language skills were let down 
only by a reticence about speaking that drove him to pick ancient Greek over 
a modern language (a good choice for a mathematician, though), what fired him 
was mathematics. He was ``highly commended'' in the 1971 National Schools 
Mathematical Contest. He took math, further math and physics at ``A'' level and 
a BA (1974)  in mathematics at Trinity College Cambridge before moving into 
the theory of knots for his PhD project under Professor Cameron Gordon, still 
at Trinity, where he was a Fellow from 1978 to 1982. He completed his thesis, 
entitled Topics in Knot Theory, in 1979, but his collaboration with Professor 
Gordon was to continue into the 1980s, leading to the development of the 
Gordon-Litherland Pairing. Rick enjoyed WBR Lickorish's contention that one 
of his and Gordon's theorems was only a ``folk theorem'' because its proof was 
found at a folk music exhibition in Helsinki.

Knot theory was Rick's wormhole into a universe denied to his baffled yet 
admiring family, who looked on at his donnish Cambridge life without 
understanding the impact of the work he was doing. Topology seemed to be for 
him an X--ray eye through which he could look beneath Nature’s spume to --- in 
Yeats’s phrase --- the ghostly paradigm of things. A ghostly paradigm that 
could be reified in plastic snap--together beads from Woolworths.

At Cambridge he met his wife--to--be Alison, a botanist, artist and fellow 
student. They married in 1977, just before Elvis Presley died. He and Alison 
would go walking in the Lake District and on the Pennine Way, with Rick 
confidently striding ahead and map reading, provisioned with Kendal mint 
cake and Mars bars, and always staying in youth hostels, of which he was a 
great supporter.

After his doctorate, Rick worked briefly at the Universities of Cambridge and 
Warwick in the UK before moving to the United States. There, he lectured for 
six months at the University of Texas at Austin and in 1983 took up a teaching 
and research position at Louisiana State University in Baton Rouge. At this 
point he and Alison separated, Alison living in Cambridge (England) and Rick 
in Baton Rouge (Louisiana), but they never lost touch. Whenever he could, 
Rick went back to the UK twice a year at Christmas and in summer to visit 
family, drop in on Alison and attend the Cambridge folk festival - he loved 
traditional music as well as some slightly left-field rock such as the Bonzo 
Dog Band and Captain Beefheart. On these visits he liked nothing better than 
to sit in a quiet, traditional pub with a pint  of real ale and the Guardian 
cryptic crossword, followed by spending the afternoon browsing the bookshops 
and the evening watching quiz shows followed by a quick late-evening pint at 
the local pub (and participating in the pub quiz).

Apart from secondments to Rice University (Houston) and the University of 
Iowa, Rick stayed at LSU until his retirement in 2022, achieving full 
professorship in 1994. His evaluations described him as a conscious and 
effective teacher and a highly respected researcher.

Rick throve on Louisiana's music, food and weather, although the state's 
right-leaning political climate suited him less well. He reveled in Mardi 
Gras and the New Orleans Jazz Festival, and here his musical tastes broadened 
to encompass zydeco, Cajun, funk and swamp rock.  He loved to retell the 
anecdotes of the Huey and Earl Long eras. He lived in a rented house that 
was the former granny flat of a bigger house in Southdowns in Baton Rouge, 
and never thought to move to a fancier place, seeming happy to settle for 
modest premises in a part of town from which he could reach his LSU office 
in a 30--minute stroll through magnolia--bright avenues and past the 
turtle--brimming University and City Park lakes. Wealth and possessions were 
unimportant to him; he spent his money on books and music and gave generously 
to charitable organizations despite being fond of quoting Clement Attlee's 
view of charity as a ``cold, grey, loveless thing''.

\[\begin{array}{c}
\includegraphics{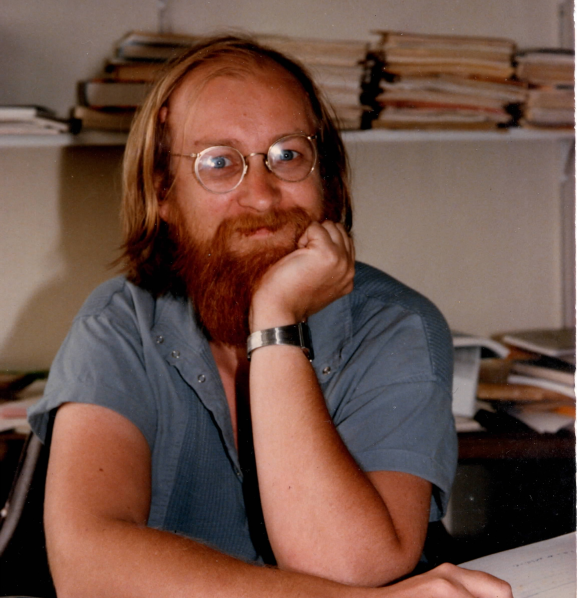} \\
\mathrm{R.\ A.\ Litherland\ in\ the\ early\ 1980s}
\end{array}\]
%[office-1980s.jpeg could go here. Caption something like: R. A. Litherland in the early 1980s]

In Baton Rouge he met his partner--to--be Cheryl Nelson, who soon moved into 
the flat with him. Rick was devastated when Cheryl died in 2008 aged only 39. 
He continued living alone in the same flat.

Rick had an interest in computing that began in adolescence with two 
educational analog quasi-computers of the 1960s, the Digi-Comp$^{\mathrm{TM}}$ 
I and the Think a Dot$^{\mathrm{TM}}$. 
His Sinclair ZX80, on which he practiced squeezing a lot into 
a very little RAM, is now in the care of the Centre for Computing History in 
Cambridge, UK. For the 1999 Interactive Fiction Competition, Rick wrote a 
TADS  game called Erehwon that got some decent reviews and came 11th in the 
competition. It featured the five Platonic solids, matter and antimatter,  
a wormhole, a Hamiltonian maze, a desert with the geometry of a Klein bottle 
-- and a touching portrait of Mr Bird, the budgerigar who shared Rick and 
Cheryl’s Baton Rouge flat and for whom Rick avoided displaying the slightest 
affection in real life. Mr Bird's role in the game is to say things like ``i 
snig the body eluctrun'' and ``hewwop lippop weeb'' (the latter a reference 
intelligible only to British baby-boomers) before being annihilated in a 
matter/antimatter clash as a sacrifice to the game logic and being replaced 
in the player's inventory by a ``burden of guilt''. Here at last was some 
output that his non--topologist friends could understand, although even here 
the in-game hints came in handy.

Rick enlivened his retirement by playing with the Turing Tumble, another 
type of toy analogue computer that can be ``programmed'' using plastic cogs 
and flipflops. At the time of his death, he was amusing himself by writing 
a simulator for this bagatelle--like machine in Python.

Rick died of complications following a stroke just as he was about to fly 
back to the UK for Christmas 2022 on his first visit in several years, 
health problems and COVID restrictions having conspired since summer 2017 
to keep him from traveling to the country he continued to regard as ``home'' 
and to which he was weighing up the idea of retiring, while -- typically 
-- showing no sign of ever making up his mind.

\section{\large\textbf{Mathematical Work}}\label{MW} %author: Sam

As of October 2023, Mathscinet lists a total of 24 publications under
Rick's name, with 674 citations in 625 publications by 531 unique citing 
authors, an impressive and admirable record for any mathematician. 
Rick's mathematical work spans several decades and touches on several
areas of low-dimensional topology. There is surely too much to cover in
much detail in a short remembrance such as this, so we will content 
ourselves to review a sampling of some of Rick's most influential results.

Rick's research started off extremely well with his joint publication 
\textit{On the signature of a link} with his advisor Cameron Gordon in
\textit{Inventiones Mathematicae}, one of the most prestigious journals 
in all of mathematics \cite{GL78}. The \textit{signature} of a knot is 
an integer-valued invariant which can be computed in a number of ways,
e.g. via skein theory using the Conway polynomial, from the infinite cyclic
cover of the knot complement, or from the linking matrix. The 
\textit{Mathscinet} reviewer Wilber Whitten said of this paper:

\begin{quotation}
The results in this paper tie together, climax, and clarify the work of 
L. Goeritz, H. F. Trotter, K. Murasugi, and others \dots on the signature 
of a tame knot or link in $S^3$. In fact the authors prove that the 
signature $\sigma(L)$ of a link $L$ is the signature $\mathrm{sig}(G)$ of 
any Goeritz matrix $G$ plus a ``corrective term'', $\mu$, which one can 
easily read from a regular projection of $L$ and which depends (with $G$) 
on the projection. 
\end{quotation}

Rick's work on the subject of the signature invariant continued in
\textit{Signatures of iterated torus knots} \cite{L792} and 
\textit{Signatures of covering links} \cite{GLM}.

One thread of Rick's research program involved \textit{Surface-links}, 
compact 2-manifolds embedded in $\mathbb{R}^4$ considered up to ambient 
isotopy, or, in other words, surfaces knotted and linked in 4-space. His
work on surface-links appeared
in papers such as \textit{Deforming twist-spun knots} \cite{L791},
\textit{The second homology of the group of a knotted surface} \cite{L81}
and \textit{Symmetries of twist-spun knots} \cite{L85}.

Another thread of Rick's research program involved knots and surfaces
embedded in 3-manifolds, a topic presaging the later rise
of virtual knot theory and relevant to big questions such as the 
slice-ribbon conjecture. Rick's work on this topic appeared in papers such as 
\textit{Surgery on knots in solid tori} \cite{L793},
\textit{Surgery on knots in solid tori. {II}} \cite{L80},
\textit{Incompressible planar surfaces in {$3$}-manifolds} \cite{GL84}
and \textit{Incompressible surfaces in branched coverings} \cite{GL841}.
Surgery descriptions of 3-manifolds and related ideas figure prominently
in some of Rick's other work as well, such as \textit{Lens spaces and 
Dehn surgery} \cite{BL} and his final work, \textit{Surgery description 
of colored knots} \cite{LW}.

Rick was always open to new interesting ideas related to knots, 
from his work on spatial graph theory, i.e., the study of knotted graphs in 
$\mathbb{R}^3$ in \textit{The Alexander module of a knotted theta-curve} 
\cite{L89} and \textit{The homology of abelian covers of knotted graphs}
\cite{L99} to his paper \textit{A generalization of the lightbulb theorem 
and {PL} {$I$}-equivalence of links} on piecewise-linear knotted $n$-spheres 
in $S^n\times S^2$ \cite{L86} to more physical invariants of knots such 
\textit{Thickness of knots} \cite{LSR}.

Lest one think Rick's work was strictly about geometric topology, his work 
included a fair amount of study of algebraic structures such as global fields
in \textit{Matching {W}itts with global fields} \cite{PSCL} and of course his
work on quandles and racks in \textit{The Betti numbers of some finite racks} 
\cite{LN} and \textit{Quadratic quandles and their link invariants} \cite{L02}.

\section{\large\textbf{Once a week for four years with R.A. Litherland}}

In this section we conclude with some thoughts and memories shared by Rick's
first and last PhD students, Sam Nelson and Steven D. Wallace.

\begin{quotation}
Rick was a great mentor and friend. Humble and soft-spoken as he was, you 
would never guess that he had made such important discoveries in his early 
career with his mentor, Cameron Gordon. He was my thesis advisor when I was 
a Ph.D candidate. Our meetings were always business, but never cold or 
impersonal. He was my teacher and I will miss him. My dearest condolences to 
his family.

\hskip 3in -- Steven Wallace
\end{quotation}

\subsection{Comments by Sam Nelson, Rick's First PhD Student}

I learned a great deal of both mathematics 
and how to be a mathematician from Rick. The reading courses we did together
(Knot Theory, PL Topology and other topics) and Rick's easy-going style of 
collaboration formed the model for my interactions with my own students and 
collaborators. I recall Rick assigning me to study a theorem in our
knot theory book \cite{Lic} by theorem number and not mentioning the title
``Gordon--Litherland Form''. These days, a quarter-century later, I am still
seeing talks at international conferences about the Gordan-Litherland form.

Rick's choice of focus for our initial research effort was a paper 
on finite type invariants, also known as Vassiliev invariants, which 
also served as my introduction to the then brand-new topic of virtual knot 
theory \cite{GPV}. I think Rick has really wanted me to focus more on the 
finite type invariants side of the paper, but I was fascinated by virtual knot
theory and pursued the (initially somewhat unfashionable) topic passionately, 
perhaps to the detriment of my job prospects.

When I found that couldn't follow the proof of a result in the paper (that 
the forbidden moves are unknotting moves), Rick suggested I should try to 
find my own proof. I tried
and succeeded, and we filed away the result for future use in my dissertation.
Then one day out of the blue I received a package from Japan, wherein
Taizo Kanenoubu had written a paper in which he cited my webpage on virtual 
knots \cite{K}. I excitedly showed the paper and citaton to Rick, who nodded and
said ``Yes, but did you read the paper?'' As it turned out, the paper was 
another proof of the same result on forbidden moves, and Rick said I had to 
write it up and submit it for publication. I did, and when it was accepted by 
JKTR it became my first publication \cite{N1}.

When Rick saw me talking with another Japanese mathematician, Seiichi Kamada, 
at an AMS conference in Lafayette, he suggested I ask him for his research.  
Seiichi sent me a few papers, one of which was a draft of a joint paper with 
Scott Carter and Masahico Saito on virtual knots and quandle homology.
My first introduction to the topic of quandle theory, the paper contained a 
conjecture about the long exact sequence in quandle homology and some 
observations about certain Alexander quandles being isomorphic; with Rick's 
guidance and help, these two became the topics of my second and third 
publications and my only joint publication with Rick \cite{LN,N2}.

My friendship with Rick was warm but always felt professional. Occasionally I
would see him at the Chimes, an excellent pub just off the LSU campus, and 
say hello, but aside from conferences I didn't see him much outside the 
department. After talks at the first AMS conference session I organized around 
2007, I was having dinner with Scott Carter and Masahico Saito at the Chimes
when I said ``Back in grad school, whenever I'd come here I'd look over at 
the bar and see Rick sitting there'' to which Scott replied ``Well, isn't that
Rick over there right now?'' to which I replied ``So it is!'' and went to 
invite him to come join us. I recall Rick did once comment about my 
choice of DJ name, ``DJ Absinthe Minded'', in my email signature, saying he 
liked the name.

After completing my PhD and moving to California I saw Rick only a few times,
only on the few occasions when I came back to Baton Rouge. I did invite
him a few times to come to Claremont and give a talk in our Topology Seminar, 
but it never quite worked out. The final time I saw 
him in person was on a visit to Baton Rouge primarily to visit my own first 
Senior Thesis student, Jose Ceniceros, who was finishing his PhD at LSU.
I was very happy that Rick attended a talk I gave via Zoom toward the end 
of the Covid-19 pandemic when receiving a local research recognition.

\subsection{Comments by Steven D. Wallace, Rick's Final PhD Student}

% author: Steven D. Wallace}

\medskip

\noindent\textit{In memory of Rick Litherland and Tim Cochran.  We miss you both.}

\medskip

%\subsection{Relationship}

My initial experience as a Ph.D. student at Louisiana State University was a 
little different than most of my colleagues who arrived to campus at the same 
time.  I had earned a thesis master's degree from Rice before coming to LSU.  
Indeed, this meant that I had a laser focused agenda for my course of study 
that my fellow classmates had to figure out over the course of their first 
several semesters.  In fact, I came to LSU to work with R.A. Litherland 
because he was Sam Nelson's advisor, and my master's thesis was in the same 
general areas of ``quandles'' and knot theory.  Luckily, Rick agreed to be my 
advisor as soon as I arrived in Baton Rouge.

%\begin{itemize}

%\item Describe my first meeting with Dr. Litherland.

%\item Talk about taking photos at graduation with all of my boisterous family.  How clearly uncomfortable he was to do it, but how he still did it.

%\item (Memorial page message) Rick was a great mentor and friend. Humble and soft-spoken as he was, you would never guess that he had made such important discoveries in his early career with his mentor, Cameron Gordon. He was my thesis advisor when I was a Ph.D candidate. Our meetings were always business, but never cold or impersonal. He was my teacher and I will miss him. My dearest condolences to his family.

%\end{itemize}

Rick was a great person, a good friend and mentor, but he was a soft-spoken, 
private person.  Our meetings were all business, but never cold or impersonal.
We met at least once a week for four years in his office.  A small room, it 
could barely hold all of his books.  We would often have to move my seat to 
access one of the lower shelves.  He had space for maybe one piece of 
loose-leaf paper on his desk between piles of what I assumed were student 
papers he was grading, notebooks, a computer, and more books.  We often stood 
side by side at the slate chalkboard opposite his desk.

%Although I am not sure he ever noticed, I made a point to not shave on the day or two before our meetings to give him the idea that I was spending all-nighters to prepare for our meetings.  His mentoring technique for me was one of support rather than driven by day-to-day assignments.

Humble as he was, you would never guess that he had made such important 
discoveries in his early career with his advisor, Cameron Gordon.  He also had 
an impact on Louisiana State University culture.  ``The Professor,'' as he was
known at \textit{The Chimes} restaurant, had his assigned stool at the bar 
that everyone in the community recognized was his.  When not doing his best 
impression of Norm from \textit{Cheers}, he could be seen walking around 
campus since he never drove a car.

My whole time at LSU, one of Rick's jobs for the Math Department was to 
organize the Topology Colloquium.  He made sure that coffee was made, and 
the speakers were lined up.  It was behind-the-scenes work that was perfectly
 suited for Rick Litherland.  The cookies accompanying the coffee always 
included Pepperidge Farm Milanos, not flashy or extravagant, not overly sweet 
or soft, but steadfast and reliable, like the man who put them out each week, 
like his mentorship style for me.

On rare occasions, we convinced him to join other faculty and graduate 
students to treat the out-of-town speaker to some Baton Rouge Cajun food.  
Even though it wasn't always at his favorite spot, I think that he enjoyed 
those moments, too.

However, he rarely joined us when some of the other faculty and the handful 
of graduate students, like myself, who were always looking for an excuse to 
eat out would treat the out-of-town speaker to some cajun food.

After about a year and a half, I had passed all of my graduate student exams, 
but had not come up with a good dissertation topic.  He said, ``maybe it's 
time for something different'' in his slightly raspy British accent.  He 
suggested that I look on the ArXiv for something that interested me.  He gave
me some search parameters and tasked me to find two or three papers to 
consider.  This was his only timed assignment, that he gave me.  By the next 
week, we were working on the surgery equivalence of colored knots, which is an 
analog of a famous theorem of Lickorish and Wallace (no relation) \cite{}. %(reference here).  
I asked, "what do you think?"  He smiled and said, "well, it is pretty 
classic knot theory, so I like it," then he chuckled.

%\subsection{Summary of our work together}

My dissertation topic involved the surgery equivalence of colored knots, which 
is an analog of a famous theorem of Lickorish and Wallace (no relation) 
\cite{Lic} and had nothing to do with quandles.  I brought the idea to Rick, 
he grinned and chuckled slightly, and, in his gruff British accent, he said, 
``well, it is pretty classic knot theory, so \dots I like it.''

Let $p\geq 3$ be an integer.  A $p-$coloring of a knot, introduced by Fox as 
an exercise in \cite{Fox}, is a surjective homomorphism, 
$\rho: \pi_1 \left(S^3 -K\right)\rightarrow D_p$, of the knot group onto the 
dihedral group of order $2p$.  We define a $p-$colored knot as the pair 
$(K,\rho)$ where $\rho$ is a $p$-coloring for the knot $K$.  In \cite{Mo1} 
two $p-$colored knots are said to be \textit{surgery equivalent} if one may 
be obtained from the other by a finite sequence of $\pm 1$ surgeries along 
unknots that preserve colorability.  Indeed, the coloring from one colored 
knot is obtained from the other by the appropriate induced homomorphisms 
arising from the sequence of surgeries.

In \cite{Mo1}, D. Moskovich classified all the surgery equivalence classes 
of $3$-colored knots using local moves on $3$-colored knot diagrams.  After 
using a small forest's worth of scratch paper trying to extend those 
techniques for other colored knots, Tim Cochran suggested that we use the 
``bordism invariants'' that he developed several years prior (see \cite{CGO}).

The result was a complete invariant for the surgery equivalence of $p$-colored 
knots (\cite{LW,Wa}).  We established a universal upper bound of $2p$ for 
the number of equivalence classes for any $p\geq 3$.  In \cite{KrMo}, Kricker 
and Moskovich sharpened the upper bound down to $p$, thereby classifying all
 $p-$colored knot surgery equivalence classes. 

When it came time to publish our results, Rick was flattered by the notion 
that he was to be the co-author.  He had asked if his contribution warranted 
that distinction, which of course it did.  Even though his accomplishments in 
knot theory were very important, his humility was always apparent.

\subsection{Knots in Washington 49.96875}

Rick's advisor, Cameron Gordon, shared a few memories of Rick at the Knots 
in Washington 49.96875 conference at George Washington University in December 
2023, which was dedicated to Rick's memory. Sam Nelson paraphrases:

``Since we were basically the same age, our relationship was more like 
colleagues than advisor and student.''

``Rick always dressed the same way, as in the picture. Once we were with a 
group going out for drinks when the bouncer decided he didn't want to let 
Rick in. `You can't come in, we have a dress code' he said, and Rick asked 
what was wrong. `No tennis shoes' said the bouncer, but Rick replied 'I'm not 
wearing tennis shoes.' `Well, no blue jeans!' he said, and Rick replied 'I'm 
wearing pants, not jeans'. `Well, no t-shirts!' the bouncer said, and Rick 
objected `I'm not wearing a t-shirt.' Realizing he couldn't win, the bouncer 
said `Ok fine, you can come in -- but you'd better not cause any trouble!' ''

``A propos of him always looking the same, he remarked to me once that when 
he was young older people would shout at him in the street `Get your hair 
cut, you bloody hippy!,' and then when he was older young skinheads would 
shout at him `Get your hair cut, you bloody hippy!' ''

\section{Open Problems}

We conclude with a few problems which we believe to be currently open related
to Rick's research and ideas.

Steven Wallace writes: 
``I was thinking about the `audioactive sequences' that Lou Kauffman 
referenced at Knots in Washington.  They were originally studied by 
Robert Morris, and then also by John Conway \cite{Conway1987} 
%(paper in 1986 called ‘The Weird and Wonderful Chemistry of Audioactive Decay’).  
Is the property that arises from the seed digit 3 unique to only that seed 
digit?  Conway found 92 `elements' for the `look-and-say sequences' but that 
may have been a different way to parse them.''

Jozef Przytycki writes:
``I have the following question related to Rick's work:
Can you find the signature of a knot using Fox coloring matrix 
(Mattman-Solis call this \textit{crossing matrix}),
that is a matrix of universal Fox coloring group (of course some 
correcting factors are allowed)?''

Sam Nelson writes: ``The Gordan-Litherland pairing has recently been extended 
to links in thickened surfaces in \cite{BCK}. How about analogues in other 
categories of knotted objects such as spatial graphs, surface-links and 
pseudoknots? In \cite{L02} Rick proposed two conjectures about quandle cocycles 
invariants of a specific form being determined by the Alexander module; in 
the same paper he proved the conjectures for torus links and 2-bridge knots.
Are they true for all knots?''

\bibliography{ricksbio}{}
\bibliographystyle{abbrv}

\bigskip

\medskip

\noindent
\textsc{Department of Mathematical Sciences \\
Claremont McKenna College \\
850 Columbia Ave. \\
Claremont, CA 91711}

\medskip

\noindent 
\textsc{Department of Mathematics \& Statistics \\
Middle Georgia State University \\
100 University Parkway \\
Macon, Georgia 31206}

\end{document}